\documentclass[12pt]{article}%
\usepackage{graphicx}
\usepackage{amsmath}
\usepackage{amsfonts}
\usepackage{amssymb}%
\providecommand{\U}[1]{\protect\rule{.1in}{.1in}}
\newtheorem{theorem}{Theorem}

\newtheorem{proposition}[theorem]{Proposition}

\begin{document}

\title{\textbf{Plinths and Pedestals}}
\author{Oleg Ogievetsky$^{\sharp,\flat,\dag}$ \& Senya Shlosman$^{\natural
,\sharp,\flat,\ddag}$\\$^{\natural}$Krichever Center for Advance Studies, Moscow, Russia;\\$^{\sharp}$Aix Marseille Univ, Universite de Toulon, \\CNRS, CPT, Marseille, France;\\$^{\flat}$Inst. of the Information Transmission Problems, \\RAS, Moscow, Russia;\\$^{\dag}$Lebedev Physical Institute, Moscow, Russia,\\$^{\ddag}$BIMSA, Beijing, China\\Oleg.Ogievetsky@gmail.com, shlosman@gmail.com}
\maketitle

\begin{abstract}
In this paper we explain what are the plinths and the pedestals of the
skyscrapers (=plane partitions), and how one can use them in order to count
the skyscrapers.

\end{abstract}

\section{Plinths}

Our skyscrapers will stand on Young diagrams, $\lambda,$ which are tables with
$\ell\left(  \lambda\right)  $ rows. The rows have descending lengths
$\lambda_{1}\geq\lambda_{2}\geq...\geq\lambda_{\ell\left(  \lambda\right)  }.$
The number of cells in $\lambda$, which equals to $\lambda_{1}+\lambda
_{2}+...+\lambda_{k},$ is denoted by $\left\vert \lambda\right\vert .$ A cell
$\left(  i,j\right)  $ of $\lambda$ is located on the $i$-th position in the
$j$-th row, so $i\leq\lambda_{j},$ and we think of $i$ increasing to the
right, while $j$ increases downward. More generally, if $\mu\subset\lambda$
are two Young diagrams, then the complement $\lambda\setminus\mu$ is called a
skew Young diagram, $\lambda/\mu.$ A skyscraper $T=\left(  T_{ij}\right)  $ is
an integer-valued function on the cells of $\lambda/\mu$ (i.e. an array of
integers of shape $\lambda/\mu$).

The array $Q=\left\{  Q_{ij},\left(  i,j\right)  \in\lambda/\mu\right\}  $ is
called a Standard Young Tableaux (SYT) if the integers $Q_{ij}$ take all the
values $\left\{  1,2,...,n=\left\vert \lambda/\mu\right\vert \right\}  $ (each
once) and are strictly increasing in every row and in every column. SYT is the
same as a linear order on the cells of $\lambda/\mu,$ consistent with the
natural partial order on $\lambda/\mu.$ We denote this finite set by
$SYT\left(  \lambda/\mu\right)  $.

The array $T=\left\{  T_{ij},\left(  i,j\right)  \in\lambda/\mu\right\}  $ of
\textbf{non-negative} integers of shape $\lambda/\mu$ is called a Semistandard
Young Tableaux (SsYT) of shape $\lambda/\mu$ (i.e., $1\leq j\leq\ell\left(
\lambda\right)  ,$ $\mu_{j}<i\leq\lambda_{j}$), if they are weakly increasing
in every row and strictly increasing in every column. We denote by $\left\vert
T\right\vert $ the volume $\sum_{\left(  i,j\right)  \in\lambda/\mu}T_{ij}$ of
the skyscraper $T,$ and by $SsYT\left(  \lambda/\mu\right)  $ the set of all
SsYT of shape $\lambda/\mu.$

The starting point of our study is the following relation, expressing the
generating function of the Semistandard Young Tableaux of shape $\lambda/\mu$
via the major index $\mathsf{maj}(\ast)$ (defined in $\left(  \ref{44}\right)
$ below) of Standard Young Tableaux of the same shape:

\begin{proposition}
(\cite{St}, Proposition 7.19.11) The generating function%
\begin{equation}
G_{\lambda/\mu}^{Ss}\left(  q\right)  \equiv\sum_{T\in SsYT\left(  \lambda
/\mu\right)  }q^{\left\vert T\right\vert }=\frac{\sum_{Q\in SYT\left(
\lambda/\mu\right)  }q^{\mathsf{maj}\left(  Q\right)  }}{\left(  1-q\right)
\left(  1-q^{2}\right)  ...\left(  1-q^{\left\vert \lambda/\mu\right\vert
}\right)  }. \label{21}%
\end{equation}

\end{proposition}

The index $\mathsf{maj}\left(  Q\right)  $ is defined as follows$.$ For every
$k=1,...,n=\left\vert \lambda/\mu\right\vert $ there is a unique cell $\left(
i,j\right)  $ for which $Q_{ij}=k.$ In such a case we will say that $\left(
i,j\right)  =Q\left(  k\right)  ,$ and we say that $k$ is the content of the
cell $\left(  i,j\right)  .$ The value $k,$ $k=1,...,n-1$ (and the cell
$Q\left(  k\right)  $) is called a descent of $Q,$ if
\begin{equation}
Q\left(  k\right)  =\left(  i,j\right)  ,\ Q\left(  k+1\right)  =\left(
i^{\prime},j^{\prime}\right)  \text{ and }j^{\prime}>j. \label{48}%
\end{equation}
In words, the entry $k$ is a descent iff the row of the entry $k+1$ is below
the row of the entry $k.$ We denote by $D\left(  Q\right)  \subset\lambda/\mu$
the set of all descents cells of $Q$. Let
\begin{equation}
\mathsf{Des}\left(  Q\right)  =\left\{  i_{1},...,i_{l}\right\}  \label{49}%
\end{equation}
be the contents of descent cells of $Q.$ By definition,%
\begin{equation}
\mathsf{maj}\left(  Q\right)  =\sum_{k:Q\left(  k\right)  \in D\left(
Q\right)  }k\equiv\sum_{k\in\mathsf{Des}\left(  Q\right)  }k. \label{44}%
\end{equation}

Looking on the relation $\left(  \ref{21}\right)  $, one recognizes in the
factor $\frac{1}{\left(  1-q\right)  \left(  1-q^{2}\right)  ...\left(
1-q^{n}\right)  }$ the well-known generating function of the set
$\mathcal{Y}_{n}$ of partitions:
\[
\mathcal{Y}_{n}=\left\{  Y=\left(  y_{1},...,y_{n}\right)  :0\leq y_{1}%
\leq...\leq y_{n}\right\}
\]
of non-negative integers into $n$ non-negative parts:%
\[
G_{\mathcal{Y}_{n}}\left(  q\right)  \equiv\sum_{Y\in\mathcal{Y}_{n}%
}q^{\left\vert Y\right\vert }=\frac{1}{\left(  1-q\right)  \left(
1-q^{2}\right)  ...\left(  1-q^{n}\right)  },
\]
where $\left\vert Y\right\vert =y_{1}+...+y_{n}.$ One would like to interpret
the remaining factor in $\left(  \ref{21}\right)  $ also as a generating
function, counting something interesting.

The remarkable fact is that there exists indeed a finite subset $Pl\left(
\lambda/\mu\right)  \subset SsYT\left(  \lambda/\mu\right)  $ of the set of
SsYT and a bijection $B^{Ss}:SsYT\left(  \lambda/\mu\right)  \rightarrow
Pl\left(  \lambda/\mu\right)  \times\mathcal{Y}_{n},$ respecting the volumes:
if $B^{Ss}\left(  T\right)  =\mathsf{p}\times Y,$ then $\left\vert
T\right\vert =\left\vert \mathsf{p}\right\vert +\left\vert Y\right\vert .$
That implies that our generating function factors:%
\[
G_{\lambda/\mu}^{Ss}\left(  q\right)  =G_{\lambda/\mu}^{Pl}\left(  q\right)
\times G_{\mathcal{Y}_{n}}\left(  q\right)  ,
\]
where of course $G_{\lambda/\mu}^{Pl}\left(  q\right)  =\sum_{\mathsf{p}\in
Pl\left(  \lambda/\mu\right)  }q^{\left\vert \mathsf{p}\right\vert }.$ The
SsYT $\mathsf{p}$-s forming $Pl\left(  \lambda/\mu\right)  $ are called
plinths. They were introduced in \cite{OS2}. Clearly, $G_{\lambda/\mu}%
^{Pl}\left(  q\right)  =\sum_{Q\in SYT\left(  \lambda/\mu\right)
}q^{\mathsf{maj}\left(  Q\right)  },$ and we get thus a bijective proof of the
Stanley formula $\left(  \ref{21}\right)  .$

The set of plinths $Pl\left(  \lambda/\mu\right)  =\left\{  \mathsf{p}%
\right\}  $ turns out to be in bijective correspondence with the set of
standard Young tableaux $SYT\left(  \lambda/\mu\right)  =\left\{  Q\right\}
.$ The plinth $\mathsf{p}\left(  Q\right)  $ is defined as follows: for the
cell $\left(  i,j\right)  =Q\left(  k\right)  $ we put%
\begin{align*}
\left[  \mathsf{p}\left(  Q\right)  \right]  _{ij}  &  =\sum_{\substack{i\in
\mathsf{Des}\left(  Q\right)  :\\i<k}}1\\
&  \equiv\text{number of descents in }Q\text{ before }k.
\end{align*}

The bijection $B^{Ss}$ works as follows: to a pair $\left\{  \mathsf{p}\left(
Q\right)  ,Y=\left(  y_{1},...,y_{n}\right)  \right\}  $ it corresponds the
SsYT $T$ via the formula%
\[
T_{ij}=\left[  \mathsf{p}\left(  Q\right)  \right]  _{ij}+y_{k},\text{ where
}\left(  i,j\right)  =Q\left(  k\right)  ,
\]
see \cite{OS2} for details.

\subsection{Mahonian statistics}

As we already said,%
\begin{equation}
G_{\lambda/\mu}^{Pl}\left(  q\right)  =\sum_{Q\in SYT\left(  \lambda
/\mu\right)  }q^{\mathsf{maj}\left(  Q\right)  }=\sum_{\mathsf{p}\in Pl\left(
\lambda/\mu\right)  }q^{\left\vert \mathsf{p}\right\vert }=\sum_{Q\in
SYT\left(  \lambda/\mu\right)  }q^{\left\vert \mathsf{p}\left(  Q\right)
\right\vert }. \label{71}%
\end{equation}
However it is not the case that the volume of the plinth $\left\vert
\mathsf{p}\left(  Q\right)  \right\vert $ equals to the index $\mathsf{maj}%
\left(  Q\right)  ,$ though the two statistics -- $\mathsf{maj}\left(
\ast\right)  $ and $\left\vert \mathsf{p}\left(  \ast\right)  \right\vert $ --
are the same (which is expressed by $\left(  \ref{71}\right)  $). So the
statistics $\left\vert \mathsf{p}\left(  \ast\right)  \right\vert $ can be
called a Mahonian statistics, in honour of Percy Alexander MacMahon. This
coincidence of statistics is demystified via the use of the Sch\"{u}tzenberger
involution $Sch:SYT\left(  \lambda/\mu\right)  \rightarrow SYT\left(
\lambda/\mu\right)  .$ It turns out that
\[
\mathsf{maj}\left(  Q\right)  =\left\vert \mathsf{p}\left(  Sch\left(
Q\right)  \right)  \right\vert .
\]

\subsection{Sch\"{u}tzenberger involution}

Initially the Sch\"{u}tzenberger involution is defined only for straight Young
tableaux $\lambda$ (i.e. when $\mu=\varnothing$). The Sch\"{u}tzenberger's
definition uses the so-called Jue-de-Taquin moves, see \cite{St}. A short
(though implicit) way of defining it is via the Robinson--Schensted (RS)
correspondence: Let $Q$ be a SYT of shape $\lambda,$ and $P$ be some other SYT
of the same shape. Let $\sigma=\left(  \sigma_{1},...,\sigma_{n}\right)  $ be
the unique permutation of $n=\left\vert \lambda\right\vert $ elements which
RS-corresponds to the pair $\left(  P,Q\right)  .$ Consider the pair $\left(
P^{\prime},Q^{\prime}\right)  $ which RS-corresponds to the reversed
permutation $\sigma^{\prime}=\left(  \sigma_{n},...,\sigma_{1}\right)  .$ Then
$Q^{\prime}=Sch\left(  Q\right)  .$

The extension of the Sch\"{u}tzenberger involution to the case of SYT of skew
shape$:Sch:SYT\left(  \lambda/\mu\right)  \rightarrow SYT\left(  \lambda
/\mu\right)  $ is given in \cite{OS2}.

\section{Pedestals}

Here we start with more general setting and we will consider general posets
instead of Young tableaux. Let $X=\left\{  \alpha_{1},...,\alpha_{n}\right\}
$ be a partially ordered set with the partial order $\preccurlyeq.$ We will
still call the $\alpha$-s the cells of $X.$ A function $T:X\rightarrow\left\{
0,1,2,...\right\}  $ is called an $X$-partition if for any pair $\alpha
_{i},\alpha_{j},$ satisfying $\alpha_{i}\preccurlyeq\alpha_{j}$ we have
$T\left(  \alpha_{i}\right)  \leq T\left(  \alpha_{j}\right)  .$ We denote by
$\mathcal{P}_{X}$ the set of all $X$-partitions. (This is a generalization of
the reversed plane partitions, \cite{St}.) We will denote by $\left\vert
T\right\vert $ the volume of $T:$
\[
\left\vert T\right\vert =\sum_{i=1}^{n}T\left(  \alpha_{i}\right)  .
\]

An $X$-partition is called standard if it is a bijection of $X$ to the set
$\left[  1,...,n\right]  .$ Clearly, if $P$ is a standard $X$-partition, then
for any pair $\alpha_{i},\alpha_{j},$ satisfying $\alpha_{i}\preccurlyeq
\alpha_{j}$ we have $P\left(  \alpha_{i}\right)  <P\left(  \alpha_{j}\right)
.$ Standard $X$-partitions are just the linear extensions of the partial order
$\preccurlyeq.$ We denote by $\mathcal{P}_{X}^{S}$ the finite set of all
standard $X$-partitions.

Finally, let $F:X\rightarrow\left[  1,...,k\right]  ,$ $k\leq n,$ be a
surjective map respecting the partial order $\preccurlyeq.$ Such a map is
called a filter. We say that $\alpha_{j}$ is below $\alpha_{i}$ iff
$\alpha_{i}\preccurlyeq\alpha_{j}\ $and $F\left(  \alpha_{i}\right)  <F\left(
\alpha_{j}\right)  .$ We say that $\alpha_{i}$ is located at the $l$-th floor
iff $F\left(  \alpha_{i}\right)  =l.$ We call an $X$-partition $T$
semistandard, if for any pair $\alpha_{i}\preccurlyeq\alpha_{j}\ $with
$F\left(  \alpha_{i}\right)  <F\left(  \alpha_{j}\right)  $ we have $T\left(
\alpha_{i}\right)  <T\left(  \alpha_{j}\right)  .$ We denote by $\mathcal{P}%
_{X,F}^{Ss}\subset\mathcal{P}_{X}$ the subset of all semistandard
$X$-partitions. (This is a generalization of the semistandard Young tableaux.)
In case the filter $F_{0}$ is a trivial one: $F_{0}\left(  \alpha_{i}\right)
=1$ for all $i,$ we clearly have $\mathcal{P}_{X,F_{0}}^{Ss}=\mathcal{P}_{X}.$

We will be interested in generating functions of the $X$-partitions and
semistandard $X$-partitions:%

\[
G_{X}\left(  q\right)  =\sum_{T\in\mathcal{P}_{X}}q^{\left\vert T\right\vert
};\ \ G_{X,F}\left(  q\right)  =\sum_{T\in\mathcal{P}_{X,F}^{Ss}}q^{\left\vert
T\right\vert }.
\]

To count the $X$-partitions we will need pedestals. Let us fix a standard
$X$-partition $P\in\mathcal{P}_{X}^{S}$ -- i.e. a linear order on $\left(
X,\preccurlyeq\right)  $, and let $Q$ be another such linear extensions of
$\preccurlyeq.$

The cell $Q^{-1}\left(  k\right)  $ is called a $P$-ascent of $Q,$ if
\[
Q^{-1}\left(  k\right)  =\alpha,\ Q^{-1}\left(  k+1\right)  =\alpha^{\prime
}\text{ and }P\left(  \alpha^{\prime}\right)  <P\left(  \alpha\right)  .
\]
In words, the cell $Q^{-1}\left(  k\right)  \in X$ is a $P$-ascent of $Q$ iff
the $P$-order of the cells $Q^{-1}\left(  k\right)  ,\ Q^{-1}\left(
k+1\right)  $ is opposite to that of $Q.$ We denote by $A_{P}\left(  Q\right)
\subset X$ the set of all ascent cells of $Q$. In case $Q^{-1}\left(
k\right)  =\alpha,$ we call the value $k$ to be the content of the cell
$\alpha.$ Let
\[
\mathsf{Asc}_{P}\left(  Q\right)  =\left\{  i_{1},...,i_{l}\right\}
\]
be the contents of $P$-ascent cells of $Q.$

We define the pedestal of $Q$ with respect to $P$ as the $X$-partition
$\mathsf{d}_{P}\left(  Q\right)  $ given by
\begin{align}
\left[  \mathsf{d}_{P}\left(  Q\right)  \right]  _{Q\left(  k\right)  }  &
=\sum_{\substack{i\in\mathsf{Asc}_{P}\left(  Q\right)  :\\i<k}}1\label{13}\\
&  \equiv\text{number of }P\text{-ascents in }Q\text{ before }k.\nonumber
\end{align}
The set of all pedestals thus defined is denoted by $\mathsf{D}_{P}\left(
X\right)  .$ They were introduced in \cite{S}.

\subsection{Pedestal polynomials}

It turns out that the generating functions
\[
G_{P}\left(  q\right)  =\sum_{Q\in\mathcal{P}_{X}^{S}}q^{\left\vert
\mathsf{d}_{P}\left(  Q\right)  \right\vert }%
\]
do not depend on $P,$ so the statistics $\left\vert \mathsf{d}_{P}\left(
\ast\right)  \right\vert $ are the same for all $P.$ The generating functions
$G_{P}\left(  q\right)  $ are called pedestal polynomials, \cite{OS1}. They
enter into the generating functions of the $X$-partitions, due to the
bijections
\[
b_{St}^{P}:\mathsf{D}_{P}\left(  X\right)  \times\mathcal{Y}_{n}%
\leftrightarrow\mathcal{P}_{X},
\]
which preserves the volumes: for all $P,Q\in\mathcal{P}_{X}^{S}$ and
$Y\in\mathcal{Y}_{n}$ we have $\left\vert b_{St}^{P}\left(  \mathsf{d}%
_{P}\left(  Q\right)  ,Y\right)  \right\vert =\left\vert \mathsf{d}_{P}\left(
Q\right)  \right\vert +\left\vert Y\right\vert .$ Therefore for any $P$ we
have%
\[
G_{X}\left(  q\right)  =\sum_{T\in\mathcal{P}_{X}}q^{\left\vert T\right\vert
}\ =\left(  \sum_{Q\in\mathcal{P}_{X}^{S}}q^{\left\vert \mathsf{d}_{P}\left(
Q\right)  \right\vert }\right)  \frac{1}{\left(  1-q\right)  ...\left(
1-q^{n}\right)  }\ .
\]
The bijections $b_{St}$ are constructed in \cite{S}. They look somewhat
similar to the bijections $B^{Ss}$ discussed in the previous section. If
$\mathsf{d}_{P}\left(  Q\right)  $ is a pedestal, $Y=\left(  y_{1}%
,...,y_{n}\right)  $ is a partition $0\leq y_{1}\leq...\leq y_{n},$ then%
\[
\left[  b_{St}\left(  \mathsf{d}_{P}\left(  Q\right)  ,Y\right)  \right]
_{\alpha}=\left[  \mathsf{d}_{P}\left(  Q\right)  \right]  _{\alpha}%
+y_{k},\text{ where }k=Q\left(  \alpha\right)  .
\]

\subsection{Pedestal matrices and the Miracle of integer eigenvalues}

Consider the pedestal matrix $M^{X},$ of the size $\left\vert \mathcal{P}%
_{X}^{S}\right\vert \times\left\vert \mathcal{P}_{X}^{S}\right\vert ,$ whose
matrix elements are given by%
\[
\left(  M^{X}\right)  _{PQ}=q^{\left\vert \mathsf{d}_{P}\left(  Q\right)
\right\vert },\ \ P,Q\in\mathcal{P}_{X}^{S}%
\]
This square matrix is stochastic, as we just said: the sum $\sum
_{Q\in\mathcal{P}_{X}^{S}}q^{\left\vert \mathsf{d}_{P}\left(  Q\right)
\right\vert }$ does not depend on $P.$ What is more, it is nondegenerate, and
the miracle happen: all its eigenvalues are polynomials in $q$ with integer
coefficients! This miracle is explained in \cite{KKOPSS}. \newline

For example, let $X$ be the Young diagram $\lambda$ with 5 cells:
\begin{tabular}
[c]{|l|l|l|}\hline
$\ast$ & $\ast$ & $\ast$\\\hline
$\ast$ & $\ast$ & \\\hline
\end{tabular}
. There are 5 standard Young tableaux $P$ of shape $\lambda,$ and the pedestal
matrix $M^{\lambda}$ is given by%

\[
M^{\lambda}=\left(
\begin{array}
[c]{ccccc}%
1 & q^{3} & q & q^{4} & q^{2}\\
q^{3} & 1 & q^{4} & q & q^{2}\\
q & q^{4} & 1 & q^{3} & q^{2}\\
q^{4} & q & q^{3} & 1 & q^{2}\\
q^{4} & q & q^{3} & q^{2} & 1
\end{array}
\right)  .
\]
Its main eigenvalue is just the row sum of the matrix elements, $1+q+q^{2}%
+q^{3}+q^{4}.$ The other four are: $\left(  1-q\right)  \left(  1+q\right)  ,$
$\left(  1-q\right)  ^{2}\left(  1+q+q^{2}\right)  ,$ $\left(  1-q\right)
\left(  1+q\right)  \left(  1+q+q^{2}\right)  ,$ and $\left(  1-q\right)
\left(  1+q\right)  \left(  1-q+q^{2}\right)  .$

\subsection{Counting semistandard $X$-partitions}

Let us go now to the semistandard $X$-partitions.

We say that $T_{1}\preccurlyeq T_{2}$ if for any $\alpha\in X$ we have
$T_{1}\left(  \alpha\right)  \preccurlyeq T_{2}\left(  \alpha\right)  .$
Clearly, there exists a minimal element $\mathsf{t}_{X,F}\in\mathcal{P}%
_{X,F}^{Ss},$ such that $\mathsf{t}_{X,F}\preccurlyeq T$ for any
$T\in\mathcal{P}_{X,F}^{Ss}.$ Then the map $\mathcal{P}_{X,F}^{Ss}%
\rightarrow\mathcal{P}_{X}$ which sends $T$ to $T-\mathsf{t}_{X,F}$ is a
bijection, and so $G_{X,F}\left(  q\right)  =q^{\left\vert \mathsf{t}%
_{X,F}\right\vert }G_{X}\left(  q\right)  ,$ and the starting term of the
series $G_{X,F}\left(  q\right)  $ is $q^{\left\vert \mathsf{t}_{X,F}%
\right\vert }.$

\subsection{More Mahonian statistics}

Since $G_{X,F}\left(  q\right)  =q^{\left\vert \mathsf{t}_{X,F}\right\vert
}G_{X}\left(  q\right)  ,$ for any $P$ we have:%
\[
G_{X,F}\left(  q\right)  =\sum_{T\in\mathcal{P}_{X,F}^{Ss}}q^{\left\vert
T\right\vert }=\left(  \sum_{Q\in\mathcal{P}_{X}^{S}}q^{\left\vert
\mathsf{d}_{P}\left(  Q\right)  \right\vert +\left\vert \mathsf{t}%
_{X,F}\right\vert }\right)  \frac{1}{\left(  1-q\right)  ...\left(
1-q^{n}\right)  }.
\]

If we apply the last relation to the case when $X$ is a skew Young tableau,
and the filter $F$ is the row filter, given by $F\left(  i,j\right)  =j,$ then
we see that the two statistics -- $\mathsf{maj}\left(  \ast\right)  $ and
$\left\vert \mathsf{p}\left(  \ast\right)  \right\vert $ -- coincide with the
pedestal statistics $\left\vert \mathsf{d}_{P}\left(  \ast\right)  \right\vert
+\left\vert \mathsf{t}_{X,F}\right\vert $, for all $P$ (though all these
functions are different). So the statistics $\left\vert \mathsf{d}_{P}\left(
\ast\right)  \right\vert +\left\vert \mathsf{t}_{X,F}\right\vert $ are also Mahonian.

\textbf{Acknowledgements.} The work of S.S. was supported by the RSF under
project 23-11-00150.

\end{document}